\documentclass[12pt]{article}

\usepackage{amsmath, amscd}
\usepackage{amssymb}

\newtheorem{prop}{Proposition}

\newtheorem{corol}[prop]{Corollary}

\newcommand{\va}{\underline{a}}

\newcommand{\vx}{\underline{x}}
\newcommand{\vy}{\underline{y}}

\newcommand{\vr}{\underline{r}}

\newcommand{\wX}{X}

\DeclareMathOperator{\tr}{tr}
\DeclareMathOperator{\proj}{proj}

\DeclareMathOperator{\sign}{sign}
\newcommand{\trm}{\tr _{m}}

\title{Some matrices with nilpotent entries, and their determinants.
}\author{Anders
Kock}
\date{}

\begin{document}
\maketitle

The present note is really a Section in a forthcoming treatise 
\cite{Kompendium} on 
differential forms in the context of Synthetic Differential Geometry 
(elaborating on \cite{SDG}, \cite{DFIC}, \cite{BM}); 
but since the methods of this Section fall entirely within elementary 
linear algebra over a commutative ring $R$, we believe that it might 
be of more general interest, and worthwile a separate publication.

The base ring $R$ over which we work is implicitly supposed to have a rich 
supply of nilpotent 
elements, in particular elements $d\in R$ with $d^{2}=0$, since 
otherwise the theory collapses to the ``theory of 0-matrices''.

For the applications  which motivated the present research, $R$ is the 
number line in a model of Synthetic Differential Geometry (SDG), but no 
assumptions in this direction are needed for what we develop here. 
The only extra assumption on $R$ that we do make, is
that ``2 is cancellable in $R$'', meaning that 
for all $x\in R$, $x+x=0$ implies 
$x=0$. This will be a standing assumption.

\section{Matrices}
We consider a commutative ring $R$. We use the word ``vector space'' as 
synonymous with ``$R$-module'', and ``linear'' means ``$R$-linear''.
A vector space is called {\em finite dimensional} if it is linearly 
isomorphic to some $R^{n}$.

We begin by describing some equationally defined subsets of $R$, 
of $R^{n}$ (=the vector space of $n$-dimensional coordinate vectors), 
and of $R^{m\cdot n}$ (=the vector space of $m\times n$-matrices 
over $R$).

The fundamental one is $D\subseteq R$,
$$D:=\{x\in R \mid x^{2}=0\}.$$

More generally, 
for $n$ a positive integer, we let $D(n)\subseteq R^{n}$ be the 
following set 
of $n$-dimensional coordinate vectors $\vx=(x_{1}, \ldots ,x_{n})$:
$$D(n):=\{(x_{1}, \ldots ,x_{n})\in R^{n}\mid  
x_{j}x_{j'}=0 \mbox{ for all } j,j' = 1, \ldots ,n\},$$
in particular ($j=j'$), 
$x_{j}^{2} =0$, so that $D(n)\subseteq D^{n}\subseteq R^{n}$. The 
inclusion $D(n)\subseteq D^{n}$ will usually be a proper inclusion, 
except for $n=1$. Note also that $D=D(1)$. Note that if $\vx$ is in $ D(m)$, 
then so is $\lambda \cdot \vx$ for any $\lambda \in R$, in particular, 
$-\vx$ is in $D(m)$ if $\vx$ is.
In general, $D(n)$ is not  stable under addition.

The notation for $D$ and $D(n)$  is the standard one of 
SDG. The following set $\widetilde{D}(m,n)$ 
was first described in \cite{SDG} \S I.16 and \S I.18, 
with the aim of constructing a combinatorial notion of differential 
$m$-form.

\medskip

The subset $\widetilde{D}(m,n)\subseteq R^{m\cdot n}$ 
is  
the following set of $m\times n$ matrices $[x_{ij}]$ ($m, n 
\geq 2$):
\begin{equation*}
\begin{split}\widetilde{D}(m,n):=\{[x_{ij}]\in & R^{m\cdot n}
\mid x_{ij} x_{i'j'}+x_{i'j}x_{ij'}=0 \\ &\mbox{ for all } i,i' = 1, 
\ldots m \mbox{ and } j,j' = 1, \ldots ,n\}.\end{split}\end{equation*}
 -- We note that the equations defining $\widetilde{D}(m,n)$ 
are row-column symmetric; equivalently, the transpose of a matrix in 
$\widetilde{D}(m,n)$ belongs to $\widetilde{D}(n,m)$.
Also clearly any $p\times q$ submatrix of a matrix in 
$\widetilde{D}(m,n)$ belongs to $\widetilde{D}(p,q)$. For if the 
defining equations
\begin{equation}x_{ij} 
x_{i'j'}+x_{i'j}x_{ij'}=0\label{defining}\end{equation}
hold for all indices $i,i',j,j'$, they hold for any subset of them.
And since each of the equations in (\ref{defining}) only involve (at 
most) four indices $i,i',j,j'$, we see that for and $m\times n$ 
matrix  to belong to $\widetilde{D}(m,n)$ it suffices that all of its 
$2\times 2 $ submatrices belong to $\widetilde{D}(2,2)$.

If $[x_{ij}]\in \widetilde{D}(m,n)$, we get in particular, by putting 
$i=i'$ in the defining equation (\ref{defining}), that for any $j,j' = 1, \ldots ,n$
$$x_{ij}x_{ij'}+x_{ij}x_{ij'}=0.$$
Since 2 is assumed cancellable in $R$, we deduce from this equation  that
$x_{ij}x_{ij'}=0$, which is to say that the $i$th row of $[x_{ij}]$ 
belongs to $D(n)$. -- Similarly, the $j$th column belongs to $D(m)$.

The equations (\ref{defining}) defining $\widetilde{D}(m,n)$ can be reformulated in 
terms of a certain bilinear map $\beta :R^{n}\times R^{n}\to 
R^{n^{2}}$, where $\beta (\vx , \vy )$ is the $n^{2}$-tuple whose $jj'$ 
entry is $x_{j}y_{j'}+x_{j'}y_{j}$. Then an  $m\times n$ matrix $X$ 
($m, n\geq 2$) is in $\widetilde{D}(m,n)$ if and only if $\beta (\vr 
_{i},\vr _{i'})=0$ for all $i,i' = 1, \ldots ,m$ ($\vr _{i}$ denoting 
the $i$th row of $X$). 
 
Note that this description is not row-column symmetric. But it has 
the advantage of making the following observation almost trivial:

\begin{prop}If an $m\times n$ matrix $X$ is in $\widetilde{D}(m,n)$, then 
the matrix $X'$ formed by adjoining to $X$ a row 
which is a linear combination of the rows of $X$, is in 
$\widetilde{D}(m+1,n)$.
\label{triv}\end{prop}

\medskip

 (There is of course a similar Proposition for columns.) Combining   this Proposition with the observation that the rows of a 
matrix in $\widetilde{D}(p,n)$ are in $D(n)$, 
 we therefore have

\begin{prop}If $X$ is a matrix in $\widetilde{D}(m,n)$, then any row 
in $X$ is in $D(n)$, and also any linear combination of rows of $X$ 
is in $D(n)$. -- Similarly for columns.
\label{triv2}\end{prop} 

We have a ``geometric'' characterization of matrices in 
$\widetilde{D}(m,n)$, which depends on the following definition. We 
say that two vectors $\vx = (x_{1} , \ldots ,v_{n})$ and $\vy = 
(y_{1}, \ldots ,y_{n})$ in $R^{n}$ are {\em neighbors} (more 
precisely, {\em first order neighbours}) if $\vx - \vy \in D(n)$.
It is clearly a reflexive and symmetric relation. To say that $\vx 
\in D(n)$ is thus equivalent to saying that $\vx $ is a neighbour of 
the zero vector $0\in R^{n}$. (This ``neigbour''-relation is closely 
related to ``the first neighbourhood of the diagonal'' known for 
schemes in algebraic geometry, see e..g.\ \cite{BM};
 this is a fundamental relation in 
SDG.)

The geometric characterization of $\widetilde{D}(m,n)$ is now the 
equivalence of 1) and 2) (or of 1) and 3)) in the following

\begin{prop}Given an $m\times n$ matrix $X=[x_{ij}]$ ($m,n \geq 2$). Then the 
following three conditions are equivalent: 1) the matrix belongs to 
$\widetilde{D}(m,n)$; 2) each of its rows is a neigbour of 
$0\in R^{n}$, and any two rows are mutual neighbours; 3)
each of its columns is a neigbour of 
$0\in R^{m}$, and any two columns are mutual neighbours.
2') any linear combination of the 
rows of $X$ is in $D(n)$; 3') any linear combination of the 
columns of $X$ is in $D(m)$. 
\label{geometric}\end{prop}

\noindent {\bf Proof.} 
We have already observed (Proposition \ref{triv2}) that 1) implies 
2'), which in turn trivially implies 2).

Conversely, assume the condition 2). Let $\vr _{i}$ denote the $i$th 
row of the matrix. Then the condition 2) in particular says 
that the $\vr _{i}$ and $\vr _{i'}$ are neighbours; 
this means that for any pair of column 
indices $ j,j'$,  
$$(\vr _{i}-\vr _{i'})_{j}\cdot (\vr _{i}-\vr _{i'})_{j'}=0$$
where for a vector $\vx \in R^{n}$, $\vx _{j}$ denotes its $j$th 
coordinate. So  $(x_{ij} 
-x_{i'j})\cdot (x_{ij'}-x_{i'j'})=0$. Multiplying out, we get  
\begin{equation}
x_{ij}x_{ij'}-x_{ij}x_{i'j'} -x_{i'j}x_{ij'} + 
x_{i'j}x_{i'j'}=0.\label{cond2}\end{equation}
The first term vanishes because $\vr _{i}\in D(n)$, and the last term 
vanishes because $\vr _{i'}\in D(n)$. The two middle terms therefore vanish 
together, proving  that the defining equations (\ref{defining}) for 
$\widetilde{D}(m,n)$ hold for the matrix.
 This 
proves equivalence of 1), 2), and 2'). The equivalence of 1),  3), 
and 3') now 
follows because of the row-column symmetry of the equations defining 
$\widetilde{D}(m,n)$.

\medskip

{\bf Remark.} The condition 2) in this Proposition was the motivation 
for the consideration of $\widetilde{D}(m,n)$, since the condition 
says that the $m$ rows of the matrix, together with the zero row, 
form an {\em infinitesimal $m$-simplex}, i.e.\ an $m+1$-tuple of
 mutual neighbour 
points, in $R^{n}$; see \cite{SDG} I.18 and \cite{DFIC}. 
(In the context of SDG, the  theory of differential $m$-forms, in its 
combinatorial formulation, 
has for its basic input-quantities such infinitesimal $m$-simplices. 
The notion of infinitesimal $m$-simplex, and of affine combinations
of the vertices of such, make invariant sense in any 
manifold $N$, due to some of the algebraic stability properties (in 
the spirit of 
Proposition \ref{xlin} below) which $\widetilde{D}(m,n)$ enjoys.)

\section{Stability properties}

We begin with a ``coordinate free'' characterization of 
$D(n)\subseteq R^{n}$. Recall that we assume that 2 is cancellable in 
$R$. (Another characterization is given in Proposition 
\ref{quadratic} below.)

\begin{prop}Let $\vx \in R^{n}$. Then $\vx \in D(n)$ if and 
only if for any linear $\alpha : R^{n}\to R$, $\alpha (x) \in D$. 
\label{dn}\end{prop}
{\bf Proof.} Assume $\vx \in D(n)$. Let $\alpha$ have matrix $(a_{1}, 
\ldots ,a_{n})$, so that $\alpha (\vx )= \sum _{j} a_{j}x_{j}$. Then 
$$(\alpha (\vx ))^{2}= (\sum _{j} a_{j}x_{j})(\sum _{j'} 
a_{j'}x_{j'}),$$
which is a sum of $n^{2}$ terms $a_{j}x_{j}a_{j'}x_{j'} = 
a_{j}a_{j'}x_{j}x_{j'}$, each of which vanish because $x_{j}x_{j'}=0$.

Conversely, assume $\alpha (\vx )\in D$ for all linear 
$\alpha :R^{n}\to R$. Taking $\alpha$ to be $\proj _{j}$ (=projection onto the $j$th 
coordinate), the assumption  gives that  $x_{j}^{2}=0$. Then taking 
$\alpha $ to be $\proj 
_{j} + \proj _{j'}$, the assumption  gives that  
$(x_{j}+x_{j'})^{2}=0$. In view of $x_{j}^{2}=0$ and $x_{j'}^{2}=0$, 
this says $2x_{j}x_{j'}=0$, and since 2 is cancellable, 
$x_{j}x_{j'}=0$.

\medskip

The following is an immediate Corollary:

\begin{prop}Let $f: R^{n}\to R^{m}$ be a linear map. Then $f$ maps 
$D(n)$ into $D(m)$.
\end{prop}
{\bf Proof.} Let $\vx \in D(n)$. To see that $f(\vx )\in D(m)$, it 
suffices, by Proposition \ref{dn}, to see that for any linear 
functional $\alpha 
:R^{m}\to R$, we have $\alpha (f(\vx ))\in D$. But $\alpha \circ f$ 
is a linear functional on $R^{n}$, and thus takes $\vx$ into $D$, by 
the Proposition \ref{dn} again.

\medskip

The set of matrices $\widetilde{D}(m,n)$ was defined for $m,n \geq 
2$ only, but it will make statements easier if we extend the 
definition by putting
$\widetilde{D}(1,n)=D(n), \widetilde{D}(m,1)=D(m)$ (here, of course, 
we identify $R^{p}$ with the set of $1\times p$ matrices, or $p\times 
1$ matrices, as appropriate). By Proposition \ref{triv2}, the 
assertion that $p\times q$ submatrices of matrices in 
$\widetilde{D}(m,n)$ are in $\widetilde{D}(p,q)$ retains its validity, 
also for $p$ or $q$ $=1$.

\begin{prop}Let $X \in \widetilde{D}(m,n)$. Then for any $p\times m$ matrix $P$, 
$P\cdot X \in \widetilde{D}(p,n)$; and for any $n\times q$-matrix $Q$, $X\cdot Q 
\in \widetilde{D}(m,q)$. 
\label{ideal}\end{prop}
{\bf Proof.} Because of the row-column symmetry of the property of 
being in $\widetilde{D}(k,l)$, it suffices to prove one of 
the two statements of 
the Proposition, say, the first. So consider the $p\times n$ matrix 
$P\cdot X$. Each of its rows is a linear combination of rows from 
$X$, hence is in $D(n)$, by Proposition \ref{triv2}.
 But also any linear combinatinon of rows in 
$P\cdot X$ is in $D(n)$, since a linear combination of linear 
combinations of some vectors is again a linear combination of these 
vectors. So the result follows from Proposition \ref{geometric}.

\medskip
Here is an alternative characterization of $D(n)\subseteq R^{n}$:
\begin{prop}Let $\vx \in R^{n}$. Then the following conditions are 
equivalent:

1) $\vx \in D(n)$;

2) for any bilinear $\phi :R^{n}\times R^{n}\to R$, $\phi (\vx ,\vx )=0$;

3) for any symmetric bilinear $\psi :R^{n}\times R^{n}\to R$, $\psi 
(\vx ,\vx )=0$.
\label{quadratic}\end{prop}
{\bf Proof.} Any bilinear $\phi :R^{n}\times R^{n}\to R$ may be 
written $\psi + \phi _{a}$ with $\psi$  bilinear symmetric and $\phi 
_{a}$ bilinear alternating, in particular, $\phi _{a}(\vy ,\vy )=0$ 
for any $\vy $. Therefore, 2) and 3) are equivalent. Assume 2). For any 
pair of indices $i,i'=1, \ldots ,n$, we have the bilinear map
\begin{equation}(\vx , \vy )\mapsto x_{i}\cdot 
y_{i'}.\label{ii}\end{equation} The assumption 2) applied to this 
bilinear map and to the given $\vx$ gives that $x_{i}\cdot x_{i'}=0$ for 
all such pairs $i,i'$, and this is the defining set of equations for 
$D(n)$, so $\vx \in D(n)$, proving 1). Finally, 1) implies 2), since 
any bilinear $R^{n}\times R^{n}\to R$ is a linear combination of the 
special bilinear maps listed in (\ref{ii}).

\section{Coordinate free aspects}

Consider an arbitrary 
vector space (= $R$-module) $V$. We let $D_{s}(V)\subseteq V$ be the set 
defined by
\begin{equation*}\{ v\in V\mid \exists \mbox{ linear } f:R^{n}\to V 
\mbox{ (for some $n$) and } \exists \vx \in D(n)\mbox{ with } f(\vx 
)=v\}.
\end{equation*}
\medskip
Also, we let $D_{w}(V)\subseteq V$ be the set defined by
\begin{equation}\{v\in V\mid \forall \mbox{ linear }\phi : V\to R, 
\quad 
\phi (v)\in D\}.
\label{Dw}\end{equation}
From Proposition \ref{dn} follows immediately that $D_{s}(V)\subseteq 
D_{w}(V)$ (whence the subscripts $s$ and $w$, for ``strong'' and 
``weak'').
 However,
\begin{prop}If $V$ is finite dimensional (i.e.\ if $V\cong R^{m}$ for 
some $m$), $D_{s}(V) = D_{w}(V)$ (denoted $D(V)$); for $V=R^{m}$, $D(V) 
= D(m)$.
\label{DV}\end{prop}
(An alternative characterization of $D(V)$, in terms of quadratic 
maps, 
may be obtained from a coordinate free version of Proposition  \ref{quadratic} 
above.)

\medskip

\noindent {\bf Proof.} Since both constructions $D_{s}(-)$ and $D_{w}(-)$ 
are preserved under linear isomorphisms, it suffices to prove the 
result for $V=R^{m}$, i.e.\ to prove $D(m)=D_{s}(R^{m})=D_{w}(R^{m})$.
Clearly $D(m) \subseteq D_{s}(R^{m})$; for, the witnessing $f$ may be 
taken to be the identity map. Also $D_{s}(R^{m})\subseteq D_{w}(R^{m})$, as 
observed for a general $V$. And finally $D_{w}(R^{m})\subseteq D(m)$ 
by Proposition \ref{dn}.

\medskip

Since $m\times n$ matrices may be identified with linear maps $R^{n}\to 
R^{m}$, we would like a characterization of 
the matrices in $\widetilde{D}(m,n)$ in 
terms of the vector space $Lin (R^{n},R^{m})$.

Let $V$ and $W$ be finite dimensional vector spaces ($V\cong 
R^{n}$, $W\cong R^{m}$, say).

\begin{prop}For a linear map $F: V\to W$, the following conditions 
are equivalent:

1) for all $v\in V$, $F(v)\in D(W)$.

2) for all $v\in V$ and all linear functionals $y:W\to R$, 
$y(F(v))\in D$.

3) (if $V= 
R^{n}$, $W= R^{m}$): $F\in \widetilde{D}(m,n)$. 
\label{DVW}\end{prop}
{\bf Proof.} The equivalence of 1) and 2) follows from Proposition 
\ref{DV}, applied to $F(v)$; 3) implies 1), by Proposition 
\ref{ideal}. Finally (assming $V=R^{n}$, $W=R^{m}$), to say that 1) holds is now equivalent to 
saying that the matrix product $F\cdot v $ is in $D(m)$ for any 
$n$-dimensional column vector $v$, or, equivalently, that any linear 
combination of the columns of $F$ is in $D(m)$. This implies by 
Proposition \ref{geometric} that $F\in \widetilde{D}(m,n)$.

\medskip
For arbitrary finite dimensional vector spaces $V$ and $W$, we may 
now define a subset $\widetilde{D}(V,W) \subseteq Lin (V,W)$ by 
saying that $F\in 
\widetilde{D}(V,W)$ if the equivalent conditions 1) and 2) in the Proposition 
hold. Then $\widetilde{D}(R^{n},R^{m}) = \widetilde{D}(m,n)$ (note 
the unfortunate interchange of the order of the arguments.)
Also, under the identification of $V$ with $Lin (R,V)$, $D(V)$ gets 
identified with $\widetilde{D}(R,V)$.

Note that if $V$ and $W$ are finite dimensional, $Lin (V,W)$ is 
finite dimensional, and so $D(Lin (V,W))\subseteq Lin (V,W)$ makes 
sense; it will in general be strictly smaller than $\widetilde{D}(V,W)$; in matrix 
terms, let $V=R^{n}, W=R^{m}$, and let $A =[a_{ij}] \in  Lin 
(V,W)$. Then to say that $A\in D(Lin (V,W)$  is to say that 
$a_{ij}a_{i'j'}=0$ for all $i,i',j,j'$, which is a strictly stronger assertion 
than (\ref{defining}) (the fact that it is strictly stronger follows 
from the description of the ``generic'' matrix in $\widetilde{D}(m,n)$ 
given at the end of the next Section.) 

\medskip
Let us finally record the ``ideal-'' properties of Proposition 
\ref{ideal} when expressed in coordinate free terms; $V,W$, as well 
as $U$, $U'$, denote finite dimensional vector spaces. 

\begin{prop}Let $F\in \widetilde{D} (V,W)$. Then for any linear maps 
$P:W\to U$ and $Q: U' \to V$,
$P\circ F \circ Q \in \widetilde{D}(U',U)$. 
\end{prop}

\section{Determinants}
We now consider square matrices, say $n\times n$. They form the 
$R$-algebra $gl(n)$;   the subset $\widetilde{D}(n,n)\subseteq gl(n)$
 satisfies the ideal 
property, Proposition \ref{ideal}, (but it is not an ideal, since it is not 
stable under addition). 
 Recall that $\wX \in \widetilde{D}(n,n)$ means that the equations 
(\ref{defining}) hold. Some of the determinant theory depends only on 
a smaller set of equations, namely on the equations
\begin{equation}x_{ij} x_{i'j'}+x_{i'j}x_{ij'}=0 \label{special} 
\end{equation}for $i\neq i'$ and $j\neq j'$. For brevity, we call a 
matrix satisfying this restricted set of equations a {\em special} 
matrix. Thus, a $2\times 2$ matrix $[x_{ij}]$ is special if
$$x_{11}x_{22}+x_{12}x_{21} =0;$$
a matrix is special iff all its $2\times 2$ submatrices are special. 
Unlike matrices in $\widetilde{D}(n,n)$ (which always are nilpotent), special 
matrices may be invertible, to wit for instance the $2\times 2$ 
matrix over ${\mathbb Q}$
$$\left[ \begin{array}{rr}
1 & -1\\
1&1
\end{array}\right].$$

Recall that the {\em trace} of an $n\times n$ matrix $\wX$ is the {\em sum} of its 
diagonal entries, $\tr (\wX )= \sum _{i}x_{ii}$. The {\em product} of the diagonal entries is 
usually not very interesting, but it will be significant here; for 
brevity, we call it the {\em multiplicative trace} of the matrix,
$$\trm (\wX ):= \prod _{i} x_{ii}.$$
\begin{prop}For special  matrices (in particular for matrices in 
$\widetilde{D}(n,n)$), multiplicative trace is a 
multilinear alternating function of the columns (or of the rows)
 of the matrix.
\end{prop}
{\bf Proof.} We do the column case. Multilinearity is clear. For the 
alternating property, it suffices to see that if we interchange two 
columns of a special matrix, then the multiplicative trace changes 
sign. For simplicity of notation, let us consider interchange of the 
two first columns of a special matrix $\wX$, with resulting matrix 
$\wX '$. Then
$$\trm (\wX ) = x_{11}x_{22}u$$
where $u$ is the product $x_{33}\cdot \ldots \cdot x_{nn}$, and
$$\trm (\wX ') = x_{12}x_{21}u,$$
with the same $u$. These two expressions differ by sign, by 
(\ref{special}), and this proves the Proposition.

\medskip

Recall the standard formla for the determinant of an $n\times n$ 
matrix $\wX$,
\begin{equation}\sum _{\sigma \in {\mathfrak S}_{n}} \sign (\sigma ) 
\prod _{i=1}^{n}x_{i\sigma (i)}.
\label{two}\end{equation}
 The product in the $\sigma$th term may be viewed as $\trm (\wX 
^{\sigma})$, where $\wX ^{\sigma}$ comes about by permuting the $n$ 
columns of $\wX$ according to $\sigma$.

Thus, we can write the standard formula for the determinant of any 
$n\times n$ matrix $\wX$ as follows:
$$\det (\wX )= \sum _{\sigma} \sign (\sigma ) \trm (\wX ^{\sigma}).$$
If $\wX$ is special, it follows from the Proposition that
$$\trm (\wX ^{\sigma} )= \sign (\sigma )\trm (\wX );$$
since $\sign (\sigma )\cdot \sign (\sigma )=1$, we have that all the 
$n!$ terms in the sum (\ref{two}) are 
equal, namely equal to $\trm (\wX )$.

So we get in particular
\begin{corol}If $\wX$ is a special $n\times n$ matrix, in particular, 
if $\wX \in \widetilde{D}(n,n)$, then we have
$$\det (\wX )= n!\; \trm (\wX ).$$
\end{corol}

\medskip

{\bf Remark.} The contention of this section is that for a matrix $X\in 
\widetilde{D}(n,n)$ ($n\geq 2$), its determinant is of interest. Clearly, over 
suitable rings $R$, there do exist non-zero matrices in 
$\widetilde{D}(n,n)$, -- take e.g.\ the $n\times n$ matrix all of whose entries are 
equal to $d\in R$, where $d\in R$ has $d^{2}=0$. This matrix, 
however, has determinant zero. Do there, for suitable $R$, exist 
$X\in \widetilde{D}(n,n)$ with non-zero determinant ? The answer is 
yes, namely one may take $R$ to be the commutative $k$-algebra containing the {\em generic} 
$X\in \widetilde{D}(n,n)$ (here, $k$ is a field of characteristic 0). 
By this, we mean the $k$-algebra
$$R:=k[X_{11},X_{12}, \ldots ,X_{nn}]/J$$
obtained from the polynomial $k$-algebra in $n^{2}$ indeterminates $X_{ij}$, 
by dividing out the ideal $J$, where $J$ is  generated by the defining equations 
(\ref{defining}) for 
$\widetilde{D}(n,n)$. In this ring $R$, the matrix $[X_{ij}]$ formed by 
the indeterminates satisfies the defining equations for being in 
$\widetilde{D}(n,n)$, by construction (in fact, it is what one would 
call the {\em generic} such matrix, for $k$-algebras); and its determinant is 
non-zero, by Theorem I.16.4 in \cite{SDG}. For instance, if $n=2$, 
the theorem quoted implies that $R$, as a vector space over $k$, is 
6-dimensional, having for its basis the (classes modulo $J$ of) the 
six polynomials
$$1, X_{11}, X_{12}, X_{21}, X_{22}, \left| \begin{array}{rr}
  X_{11}& X_{12}\\ X_{21}& X_{22}\end{array}\right|.$$
More generally, the $k$-algebra $R$ containing the generic matrix $X$ in 
$\widetilde{D}(m,n)$ is finite dimensional, having for its basis the 
determinants of all $p\times p$-submatrices of $X$ (the $0\times 
0$-matrix is taken to be he constant polynomial 1); see loc.cit.\

\section{Non-linear aspects}
Assume that $g:R^{m}\to R^{l}$ is a map, not necessarily linear. Then 
if $X$ is an $m\times n$ matrix, we get an $l\times n$ matrix $g\cdot 
X$ by 
applying $g$ to each of the $n$ columns of $X$. If $g$ is linear, so 
given by an $l\times m$ matrix, $g\cdot X$ is just the standard 
matrix product of $g$ and $X$.

If $\va \in R^{n}$ (viewed as a column matrix), $X\cdot \va \in R^{m}$ is a 
linear combination of the columns of $X$ (with coefficients the 
entries of $\va$). Any linear map $g: R^{m}\to R^{l}$ preserves 
linear combinations, which in matrix theoretic formulation says
\begin{equation}g\cdot (X\cdot \va )= (g\cdot X)\cdot \va 
,\label{ass}\end{equation}
which is just the associative law for matrix multiplication. A 
crucial property of matrices $X\in \widetilde{D}(m,n)$ is the 
following Proposition:
\begin{prop}Let $X\in \widetilde{D}(m,n)$, and let $g: R^{m}\to 
R^{l}$ be a 0-preserving polynomial map. Then $g$ preserves linear 
combinations of the columns of $X$, i.e.\ the law (\ref{ass}) holds. 
\label{xlin}\end{prop}
{\bf Proof.} It is enough to consider the case where $l=1$. To say 
that $g$ is a 0-preserving polynomial map is to say that
$$g(u) = g_{1}(u)+g_{2}(u,u)+ \ldots +g_{p}(u, \ldots ,u)$$
with $g_{k}:R^{m}\times \ldots \times R^{m} \to R$ $k$-linear 
symmetric. We shall do the case of ``degree-2'' polynomials only, so
$$g(u)= g_{1}(u)+g_{2}(u,u)$$
with $g_{1}$ linear and $g_{2}$ bilinear symmetric.
Since (\ref{ass}) holds for $g=g_{1}$, it suffices to see that
it holds for the $g$ given by $u\mapsto g_{2}(u,u)$; it does so, because both sides 
of (\ref{ass}) then give 0, as we shall argue. First $X\cdot \va \in D(m)$, by 
Proposition \ref{triv2}, and it is therefore killed by $u\mapsto 
g_{2}(u,u)$, by Proposition \ref{quadratic}. On the other hand, the 
matrix $g_{2}\cdot X$ has for its columns $g_{2}(c_{j},c_{j})$, and 
since $g_{2}(-,-)$ is symmetric bilinear, these columns are all 0, 
again by Propositions \ref{triv2} and \ref{quadratic}. 

\medskip

\noindent {\bf 
Remark.} Consider for a moment the real numbers ${\mathbb R}$. 
If $g:{\mathbb R}^{m}\to {\mathbb R}^{l}$ is a smooth zero 
preserving map, then it may be written $g_{l}+h$ with $g_{l}$ linear, and
$h$ a  remainder of the form $u\mapsto g_{2}(u,u)\cdot k(u)$ with 
$g_{2}$ bilinear symmetric (and $k$ smooth). This assumption on $g$ (except the 
smoothness), makes sense also for a general commutative ring $R$ 
instead of ${\mathbb R}$. Inspecting the proof of 
Proposition \ref{xlin}, we see that we might as well have proved the 
following Proposition; we did not present it as our ``primary'' 
formulation, because its seems like a more ad hoc result. It is, 
however, in this form that it is applied in SDG. (In fact, in SDG, 
the decomposition assumed in the Proposition obtains for {\em any} zero-preserving map 
$g: R^{m}\to R^{l}$.)

\begin{prop}Let $g:R^{m}\to R^{l}$ be a zero preserving map, and 
assume $g$ may be written $g_{l}+h$ with $g_{l}$ linear, and
$h$ a  remainder of the form $u\mapsto g_{2}(u,u)\cdot k(u)$ with 
$g_{2}$ bilinear symmetric. If $X\in \widetilde{D}(m,n)$, $g$ preserves linear 
combinations of the columns of $X$, i.e.\ the law (\ref{ass}) holds.
\end{prop}

\end{document}